\documentstyle[12pt,amsfonts,amssymb,leqno]{article}

\newtheorem{proposition}{Proposition} 
\newtheorem{thm}[proposition]{Theorem}
\newtheorem{cor}[proposition]{Corollary}
\newtheorem{lemma}[proposition]{Lemma}
\newtheorem{defn}[proposition]{Definition}

\newcommand{\forces}{\ \ |\!\!\!|\!- }

\begin{document}

\begin{center}
{\Large \bf $\nabla_\kappa$, remarkable cardinals, and $0^\#$}
\end{center}

\begin{center}
\renewcommand{\thefootnote}{\fnsymbol{footnote}}
\renewcommand{\thefootnote}{arabic{footnote}}
\renewcommand{\thefootnote}{\fnsymbol{footnote}}
{\large Ralf-Dieter Schindler}${}^{a}$\footnote[1]{
\noindent
1991 {\it Mathematics Subject Classification.} Primary 03E55, 03E15. Secondary
03E35, 03E60.\\
Keywords: set theory/proper forcing/large cardinals.}
\renewcommand{\thefootnote}{arabic{footnote}}
\end{center}
\begin{center} 
{\footnotesize
${}^a${\it Institut f\"ur formale Logik, Universit\"at Wien, 1090 Wien, Austria}} 
\end{center}

\begin{center}
{\tt rds@logic.univie.ac.at}

{\tt http://www.logic.univie.ac.at/${}^\sim$rds/}\\
\end{center}

\bigskip
{\footnotesize {\bf Abstract.}
We generalize $\nabla(A)$, which was introduced in \cite{jsl}, to larger cardinals.
For a regular cardinal $\kappa > \aleph_0$ we denote by
$\nabla_\kappa(A)$ the statement that $A \subset \kappa$ and for all regular $\theta
> \kappa$ do we have that
$$\{ X \in [L_\theta[A]]^{<\kappa} \ \colon \ X \cap \kappa \in \kappa \wedge
{\rm otp}(X \cap {\rm OR}) \in {\rm Card}^{L[A \cap X \cap \kappa]} \}$$ is stationary
in $[L_\theta[A]]^{<\kappa}$.

It was shown in \cite{jsl} that $\nabla_{\aleph_1}(A)$ can hold in a set-generic
extension of $L$.
We here prove that $\nabla_{\aleph_2}(A)$ can hold in a semi-proper
set-generic extension of $L$, whereas $\nabla_{\aleph_3}(\emptyset)$ is equivalent
with the
existence of $0^\#$.}

\bigskip
Let $A \subset \omega_1$. In \cite{jsl} we introduced the following assertion, denoted
by $\nabla(A)$:
$$\{ X \in [L_{\omega_2}[A]]^\omega \colon \exists \alpha < \beta \in
Card^{L[A \cap \alpha]} \ \exists \pi \ \pi \colon L_\beta[A \cap \alpha] \cong 
X \prec
L_{\omega_2}[A] \}$$ is stationary in $[L_{\omega_2}[A]]^\omega$.
The present note is concerned with generalizations of $\nabla(A)$ to larger
cardinals. 

\begin{defn}\label{def-nabla}
Let $\kappa$ and $\theta$ both be regular cardinals, $\aleph_0 < 
\kappa < \theta$. Then by
$\nabla_\kappa^\theta(A)$ we denote the statement that $A \subset \kappa$ and
$$\{ X \in [L_\theta[A]]^{<\kappa} \ \colon \ X \cap \kappa \in \kappa \wedge
{\rm otp}(X \cap {\rm OR}) \in {\rm Card}^{L[A \cap X \cap \kappa]} \}$$ is stationary
in $[L_\theta[A]]^{<\kappa}$. By $\nabla_\kappa(A)$ we denote the statement that 
$\nabla_\kappa^\theta(A)$ holds for all regular $\theta > \kappa$.  

Moreover, we write $\nabla_\kappa^\theta$ for $\nabla_\kappa^\theta(\emptyset)$,
and 
$\nabla_\kappa$ for 
$\nabla_\kappa(\emptyset)$.
\end{defn}

It is clear that $\nabla(A)$ is $\nabla_{\aleph_1}^{\aleph_2}(A)$.
The following theorem is established by the proofs in \cite{jsl}. 

\begin{thm}\label{reptition}
Equiconsistent are:

(1) $ZFC +$ ``$L({\mathbb R})$ is absolute for proper forcings,''

(2) $ZFC + V=L[A] + \nabla_{\aleph_1}(A)$, and

(3) $ZFC +$ ``there is a remarkable cardinal.''
\end{thm}

Let us repeat the definition of a remarkable cardinal for the convenience of the
reader.

\begin{defn}\label{def-remarkable}
A cardinal $\kappa$ is called {\em remarkable} iff for all regular cardinals $\theta >
\kappa$ there are $\pi$, $M$, ${\bar \kappa}$, $\sigma$, $N$, and ${\bar
\theta}$ such that the following hold:

$\bullet$ $\pi \colon M \rightarrow H_\theta$ is an elementary embedding,

$\bullet$ $M$ is
countable and transitive,

$\bullet$ $\pi({\bar \kappa}) = \kappa$,

$\bullet$ $\sigma \colon M \rightarrow N$ is an elementary embedding with 
critical point ${\bar \kappa}$,

$\bullet$  
$N$ is
countable and transitive, 

$\bullet$ ${\bar \theta} = M \cap OR$ is a regular cardinal in $N$,  
$\sigma({\bar \kappa}) > {\bar \theta}$, and

$\bullet$ $M = H^N_{\bar \theta}$, i.e., $M \in N$ and $N \models$ 
``$M$ is the set of all sets
which are hereditarily smaller than ${\bar \theta}$." 
\end{defn}

Lemma 1.6 of \cite{jsl} gave an important characterization of remarkable cardinals.

\begin{defn}\label{def-char} {\em (\cite{jsl} Definition 1.5)}
Let $\kappa$ be a cardinal. Let $G$ be
$Col(\omega,<\kappa)$-generic over $V$, let $\theta > \kappa$ be a
regular cardinal, and let
$X \in [H_\theta^{V[G]}]^\omega$. We say that $X$ {\em condenses remarkably} if 
$X = ran(\pi)$ for some 
elementary $$\pi \colon (H_\beta^{V[G \cap
H_\alpha^V]};\in,H_\beta^V,G \cap
H_\alpha^V) \rightarrow (H_\theta^{V[G]};\in,H_\theta^V,G)$$ where 
$\alpha = crit(\pi) < \beta < \kappa$ and $\beta$ is a regular cardinal (in $V$).
\end{defn}

\begin{lemma}\label{char-remarkable} {\em (\cite{jsl} Lemma 1.6)}
A cardinal $\kappa$ is remarkable if and only if for all regular cardinals $\theta >
\kappa$ do
we have that
$$\forces_{Col(\omega,<\kappa)}^V \ {\rm `` } 
\{ X \in [H_{\check \theta}^{V[{\dot G}]}]^\omega \colon
X {\rm \ condenses \ remarkably } \} {\rm \ is \ stationary." }$$
\end{lemma}

Here is a sufficient criterion for being remarkable in $L$:

\begin{lemma}\label{char-rem}
Let $\kappa$ be a regular cardinal, and suppose that
$\nabla_\kappa$ holds.
Then $\kappa$ is remarkable in $L$.
\end{lemma}

{\sc Proof.} It is easy to see that
$\nabla_\kappa^{\kappa^+}$ implies that $\kappa$ is an
inaccessible cardinal of $L$. 

Fix $\theta > \kappa$, a regular cardinal.
By $\nabla_\kappa^{\theta^+}$, we may pick some $\pi \colon L_\gamma \rightarrow
L_{\theta^+}$ such that $\gamma < \kappa$ is a (regular) cardinal in $L$.
Let $\pi(\alpha) = \kappa$ and $\pi(\beta) = \theta$.
Let ${\bar G}$ be $Col(\omega,<\alpha)$-generic over $V$
and let
$G \supset {\bar G}$ be $Col(\omega,<\kappa)$-generic over $V$.
Then $\pi$ extends, in
$V[G]$, to some $${\tilde \pi} \ \colon \ L_\gamma[{\bar G}] \rightarrow
L_{\theta^+}[G].$$ Let ${\frak M} \in L_\gamma[{\bar G}]$ be a model of finite type
with universe $L_\beta[{\bar G}]$. We have that
$${\tilde \pi} \upharpoonright L_\beta[{\bar G}] \ \colon \ 
{\frak M} \rightarrow {\tilde
\pi}({\frak M}).$$
Notice that $\gamma < \kappa$, and therefore $L_\beta[{\bar G}]$ is countable in
$L[G]$. By absoluteness (cf. \cite{jsl} Lemma 0.2), 
there is hence some 
$\sigma \in L_{\theta^+}[G]$ such that
$$\sigma \ \colon \ {\frak M} \rightarrow {\tilde \pi}({\frak M}).$$
Therefore, $\forces^{L_{\theta^+}}_{Col(\omega,<\kappa)}$ ``there is some 
countable $X \prec 
{\tilde \pi}({\frak M})$ such that $X \cap \kappa \in \kappa$ and ${\rm otp}(X \cap
{\rm OR})$ is a cardinal in $L[{\dot G} \cap L_{X \cap \kappa}]$.'' 

Pulling this assertion
back via ${\tilde \sigma}$ yields that
$\forces^{L_{\gamma}}_{Col(\omega,<\alpha)}$ ``there is some 
countable $X \prec 
{\frak M}$ such that $X \cap \alpha \in \alpha$ and ${\rm otp}(X \cap
{\rm OR})$ is a cardinal in $L[{\dot G} \cap L_{X \cap \alpha}]$.'' As ${\frak M}$ was
arbitrary, we thus have
$\forces^{L_{\gamma}}_{Col(\omega,<\alpha)}$ ``the set of all 
$X \in [L_\beta[{\dot G}]]^\omega$ 
such that $X$ condenses remarkably is stationary.'' 
Lifting this up via
$\pi$ yields
$\forces^{L}_{Col(\omega,<\kappa)}$ ``the set of all $X \in [L_\theta[{\dot
G}]]^\omega$ 
such that $X$ condenses ramarkably is stationary.'' 

We have shown that $\kappa$ is remarkable in $L$, using Lemma \ref{char-remarkable}. 

\bigskip
\hfill $\square$ (\ref{char-rem})

\bigskip
It is easy to see that for no $\kappa$ can $\nabla_\kappa^{\kappa^+}$ hold in $L$.
We shall now consider the task of forcing $\nabla_\kappa(A)$ to hold in a
(set-) generic extension of $L$.

As to $\nabla_{\aleph_1}(A)$, ${\sf Con}(3) \Rightarrow {\sf Con}(2)$ in
Theorem
\ref{reptition} is shown by proving that if $\kappa$ is remarkable in $L$ and $G
\subset \kappa$ is $Col(\omega,<\kappa)$-generic over $L$ then $\nabla_\kappa(G)$
holds in $L[G]$. Let us now turn towards $\nabla_{\aleph_2}(A)$.

\begin{thm}\label{omega_2}
Let $\kappa$ be remarkable in $L$, and suppose that there is no $\lambda < \kappa$
such that $L_\kappa \models$ ``$\lambda$ is remarkable.'' 
There is then a semi-proper forcing ${\mathbb P} \in L$
with the property that in $V^{\mathbb P}$ 
there is some $A$ such that $\nabla_{\aleph_2}(A)$ holds. 
\end{thm}

{\sc Proof.} Let $Nm$ denote Namba forcing.
Let $\theta > \omega_2$ be regular. By
${\mathbb P}_\theta$ we shall denote the forcing
$$Col(\omega_2,\theta) \star Nm.$$ Notice that $Col(\omega_2,\theta)$ turns the
cofinality of each cardinal $\xi \in [\omega_2,\theta]$ with former cofinality $\geq
\omega_2$ into $\omega_2$, and therefore ${\mathbb P}_\theta$ turns the
cofinality of each such cardinal into $\omega$. Moreover, ${\mathbb P}_\theta$ is
semi-proper by \cite{shelah}.

We shall now define an $RCS$ iteration $({\mathbb Q}_i \colon i \leq \kappa)$ as
follows. We let ${\mathbb Q}_0 = \emptyset$, and for limit ordinals $\lambda \leq
\kappa$ we let ${\mathbb Q}_\lambda$ be the revised limit of the ${\mathbb Q}_i$,
$i<\lambda$. Now suppose that ${\mathbb Q}_i$ has been defined for some $i<\kappa$.
It will be easy to verify that inductively, $\forces^{{\mathbb Q}_i}_L$ ``${\dot G}
\subset \omega_2$.'' By Lemma \ref{char-rem} and our assumption that no $\lambda <
\kappa$ is remarkable in $L_\kappa$, 
for each $p \in {\mathbb Q}_i$ there is some (least) $\theta_p < \kappa$
such that $$\lnot \ ( p \forces^{{\mathbb Q}_i}_L \ {\rm `` }
\nabla_{\aleph_2}^{\theta_p}({\dot G}) {\rm \ holds" } ).$$
Letting
$\theta = sup_{p \in {\mathbb Q}_i} \ \theta_p <
\kappa$, we therefore have that $$\forces^{{\mathbb Q}_i}_L \ {\rm `` }
\nabla_{\aleph_2}^\theta({\dot G}) {\rm \ fails." }$$ We then let
${\dot {\mathbb Q}}$ be a name for ${\mathbb P}_{\theta}$, as being defined in 
$L^{{\mathbb Q}_i}$, and we set $${\mathbb Q}_{i+1} = {\mathbb Q}_i \star 
{\dot {\mathbb Q}}.$$

Now set ${\mathbb P} = {\mathbb Q}_\kappa$. 
Let $G$ be ${\mathbb P}$-generic over $L$. We have $G \subset
\kappa$. As ${\mathbb P}$ is
semi-proper, $\omega_1^{L[G]} = \omega_1^L$. It is moreover easy to see that $\kappa =
\omega_2^{L[G]}$. We are left with having to verify that $\nabla_\kappa(G)$
holds in $L[G]$.

Suppose not. In fact suppose that there are $p \in {\mathbb P}$ and some
(least) $\theta$ such that
$$\lnot \ ( p \forces^{\mathbb P}_L \ {\rm `` }
\nabla_{\aleph_2}^\theta({\dot G}) {\rm \ holds" } ).$$
In $L$, 
we may pick $\pi \colon L_\gamma 
\rightarrow
L_{\theta^+}$ and $\sigma \colon L_\gamma \rightarrow L_{\tilde \gamma}$ such
that $\gamma < {\tilde \gamma} < \omega_1$, $\kappa \in ran(\pi)$, $\alpha
= \pi^{-1}(\kappa)$ is the critical point of $\sigma$, $\sigma(\alpha)
> \gamma$ and $\gamma$ is a regular cardinal in $L_{\tilde \gamma}$. 

Let ${\bar {\mathbb P}} = \pi^{-1}({\mathbb P})$ and ${\tilde {\mathbb P}} =
\sigma({\bar {\mathbb P}})$. It is easy to see that ${\bar {\mathbb P}}
= {\tilde {\mathbb P}} \upharpoonright \alpha$ (with the obvious meaning).
Let $\beta = \pi^{-1}(\theta)$. Notice that, using $\pi$,
there is some $q \in {\bar {\mathbb P}}$ such that $\beta$ is least with
$$\lnot \ ( q \forces^{\bar
{\mathbb P}}_{L_\gamma} \ {\rm `` }
\nabla_{\aleph_2}^\beta({\dot G}) {\rm \ holds" } ).$$   
Therefore, there is some $\beta^\star \geq \beta$ such that forcing with
${\mathbb P}_{\beta^\star}$, as defined in $L_{\tilde \gamma}^{{\bar {\mathbb P}}}$, 
is the next step right after
forcing with ${\bar {\mathbb P}}$ in the iteration ${\tilde {\mathbb P}}$. 

Let 
${\bar G} \in L$ be ${\bar {\mathbb P}}$-generic over 
$L_\gamma$ (and hence over $L_{\tilde \gamma}$, too) such that $q \in {\bar G}$, 
and let ${\tilde G}
\supset
{\bar G}$ be ${\tilde {\mathbb P}}$-generic over $L_{\tilde \gamma}$. Then $\sigma$
lifts to ${\tilde \sigma} \colon L_\gamma[{\bar G}] \rightarrow 
L_{\tilde \gamma}[{\tilde G}]$. In order to derive a contradiction it now suffices to
prove that $\nabla^\beta_\alpha({\bar G})$ holds in $L_\gamma[{\bar G}]$.

Let ${\frak M} \in L_\gamma[{\bar G}]$ be a model 
of finite type with universe $L_\beta[{\bar G}]$. We have
$${\tilde \sigma} \upharpoonright L_\beta[{\bar G}] \ \colon \ {\frak M} \rightarrow
{\tilde \sigma}({\frak M}).$$ We would now like to build a tree $T \in 
L_{\tilde \gamma}[{\tilde G}]$ searching for an embedding like this one.

\bigskip
{\bf Claim 1.} In $L_{\tilde \gamma}[{\tilde G}]$, $L_\beta[{\bar G}] 
= \bigcup_{n<\omega} \
X_n$, where for each $n<\omega$, $X_n \subset X_{n+1}$, $X_n \in L_\gamma[{\bar G}]$, 
and ${\rm Card}(X_n) = \alpha$ in $L_\gamma[{\bar G}]$. 

\bigskip
{\sc Proof.} Let $F \colon \alpha \rightarrow \beta$, $F \in 
L_{\tilde \gamma}[{\tilde G}]$, be surjective, and let $f \colon \omega \rightarrow
\alpha$, $f \in L_{\tilde \gamma}[{\tilde G}]$, be cofinal, where $F$, $f$ are the
objects adjoined by forcing with ${\mathbb P}_{\beta^\star}$, as defined in
$L_{\tilde \gamma}^{{\bar {\mathbb P}}}$. Let $$X_n' = F {\rm " } f(n) {\rm , \ for \ }
n<\omega.$$ Notice that $F \upharpoonright \xi \in L_{\tilde \gamma}[{\bar G}]$ (and
hence $\in L_\gamma[{\bar G}]$) for each $\xi < \alpha$. In particular,
$X_n' \in L_\gamma[{\bar G}]$ for each $\xi < \alpha$. The rest is easy.  

\bigskip
\hfill $\square$ (Claim 1)

\bigskip
Now fix $(X_n \colon n<\omega)$ as provided by Claim 1. We may and shall assume that
$(X_n;...) \prec {\frak M}$ for all $n<\omega$.

\bigskip
{\bf Claim 2.} ${\tilde \sigma} \upharpoonright X_n \in L_{\tilde \gamma}[{\tilde G}]$
for each $n<\omega$.

\bigskip
{\sc Proof.} Let $f \colon \alpha \rightarrow X_n$ be bijective, $f \in L_\gamma$.
For $x \in X_n$ we'll then have that $y = {\tilde \sigma}(x)$ iff there is
some $\xi < \alpha$ with $x = f(\xi) \wedge y = {\tilde \sigma}(f)(\xi)$. But $f$
and ${\tilde \sigma}(f)$ are both in $L_{\tilde \gamma}[{\tilde G}]$. Therefore,
${\tilde \sigma} \upharpoonright X_n \in L_{\tilde \gamma}[{\tilde G}]$.  

\bigskip
\hfill $\square$ (Claim 2)

\bigskip
Now let $T$ be the tree of height $\omega$ consisting of all $(X_n,\tau)$, where 
$n<\omega$ and $\tau \colon (X_n;...) \rightarrow {\tilde \sigma}({\frak M})$ is
elementary, ordered by $(X_n,\tau) \leq (X_m,\tau')$ if and only if $n \geq m$ and
$\tau \supset \tau'$. Of course, $T \in L_{\tilde \gamma}[{\tilde G}]$.
Claim 2 witnesses that $T$ is illfounded in $V$. $T$ is hence illfounded in
$L_{\tilde \gamma}[{\tilde G}]$ as well. This buys us that in 
$L_{\tilde \gamma}[{\tilde G}]$, there is some $$\tau \ \colon \ {\frak M} \rightarrow
{\tilde \sigma}({\frak M}).$$ We thus have that
$L_{\tilde \gamma}[{\tilde G}] \models$ ``there is some $X \prec {\tilde
\sigma}({\frak M})$ such that ${\rm Card}(X) < \sigma(\alpha)$, $X \cap 
\sigma(\alpha) \in \sigma(\alpha)$, and ${\rm otp}(X \cap OR) \in
{\rm Card}^{L[{\tilde G} \cap X \cap \sigma(\alpha)]}$.'' Pulling this back via
${\tilde \sigma}$ gives that
${L_\gamma}[{\bar G}] \models$ ``there is some $X \prec 
{\frak M}$ such that ${\rm Card}(X) < \alpha$, $X \cap 
\alpha \in \alpha$, and ${\rm otp}(X \cap OR) \in
{\rm Card}^{L[{\dot G} \cap X \cap \alpha]}$.'' 

As ${\frak M}$ was arbitrary,
this shows that $\nabla^\beta_\alpha({\bar G})$ holds in ${L_\gamma}[{\bar G}]$. 
 
\bigskip
\hfill $\square$ (\ref{omega_2})

\bigskip
Our Theorem \ref{omega_2} strengthens a result which is proved in
Chapter 7 of \cite{thoralf} and which (in the terminology provided by
Definition \ref{def-nabla}) 
shows that, if $0^\#$ exists then there is a semi-proper (set-) 
forcing extension $V$ of $L$ in
which there is some $A \subset \omega_2$ such that $V = L[A]$ and
$\nabla_{\aleph_2}^{\aleph_{\omega+1}}$ hold in $V$. 

We get the following corollary to Lemma \ref{char-rem} and Theorem \ref{omega_2}.

\begin{cor}\label{corollary}
Equiconsistent are:

(1) $ZFC + V=L[A] + \nabla_{\aleph_2}(A)$, and

(2) $ZFC +$ ``there is a remarkable cardinal.''
\end{cor}

We finally turn towards $\nabla_\kappa$ for $\kappa \geq {\aleph_3}$.

\begin{lemma}\label{zero-sharp} 
Let $\kappa$ be a regular cardinal, $\kappa \geq \aleph_3$. 
Suppose that $\nabla_\kappa^{\kappa^+}$ holds. Then $0^\#$ exists.
\end{lemma} 

{\sc Proof.} Suppose not. Pick $\pi \colon L_\beta \rightarrow L_{\omega_4}$ such
that $\omega_2 < \alpha = c.p.(\pi) < \omega_3$ and $\beta$ is a cardinal of $L$.
We have that ${\cal P}(\alpha) \cap L \subset L_\beta$, and we may hence define
the ultrapower $Ult(L;U)$, where $X \in U$ iff $X \in {\cal P}(\alpha) \cap L \wedge
\alpha \in \pi(X)$. As $0^\#$ does not exist, $cf^V(\alpha^{+L}) > \omega$ as a
consequence of Jensen's Covering Lemma for $L$. By standard methods this implies
that $Ult(L;U)$ is well-founded. So $0^\#$ does exist after all. Contradiction!

\bigskip
\hfill $\square$ (\ref{zero-sharp})

\begin{lemma}\label{zero-sharp2}
Suppose that $0^\#$ exists. Then $\nabla_\kappa$ holds for every regular cardinal
$\kappa > \aleph_0$.
\end{lemma}

{\sc Proof.} We consider $0^\#$ as a subset of $\omega$.
Fix $\kappa$. Let ${\frak M} = (L_\kappa;{\vec P})$ 
be a model of finite type with universe
$L_\kappa$. Let $E \subset \kappa \times \kappa$ be well-founded and 
such that $(\kappa;E)$ condenses to $(L_\kappa[0^\#];\in)$ via the isomorphism $\sigma
\colon \kappa \rightarrow L_\kappa[0^\#]$, and let $\Sigma = 
\sigma \upharpoonright
\sigma^{-1} {\rm " } \kappa$. Let $\theta > \kappa$ be regular, and let
$$\pi \colon (L_\beta;\in,0^\#,{\bar E},{\bar \Sigma},{\vec G}) \rightarrow
(L_\theta;\in,0^\#,E,\Sigma,{\vec F})$$ be such that $\beta < \kappa$ and 
$ran(\pi) \cap \kappa \in
\kappa$. It is then straightforward to check that
for all $\gamma < \beta$, $\gamma^{+L} < \beta$. Therefore, $\beta \in {\rm Card}^L$.
As ${\frak M} = (L_\kappa;{\vec P})$ was arbitrary, this means
that $\nabla_\kappa^\theta$ holds.

\bigskip
\hfill $\square$ (\ref{zero-sharp2})

\begin{cor}
Let $\kappa \geq \aleph_3$ be a regular cardinal. Equivalent are:

(1) $\nabla_\kappa$ holds, and

(2) $0^\#$ exists. 
\end{cor}

We conclude with a few remarks. Suppose that $\kappa \geq \aleph_3$
is a regular cardinal. It can be
shown that $V=L[A] \wedge \nabla_\kappa^{\kappa^+}(A)$ implies 
that every element of $H_\kappa$
has a sharp (but of course, $A^\#$ doesn't exist
in $L[A]$). Moreover, if $A \subset \kappa$ is such that
$H_\kappa^{L[A^\#]} = H_\kappa^{L[A]}$ then $L[A] \models \nabla_\kappa(A)$. In
particular, if $L^\# = L[E]$ ($=$ the least weasel which is closed under sharps) then
$L^\# \models$ ``for all regular cardinals $\kappa > \aleph_0$, 
$\nabla_\kappa(E \cap \kappa)$ holds.''

\end{document}